\documentclass[12, reqno]{amsart}
\usepackage{amsmath}
\usepackage{amscd,amsthm,amsfonts,amsopn,amssymb,verbatim}
\parskip =\medskipamount
\numberwithin{equation}{section}
\newtheorem{theorem}{Theorem} 

\theoremstyle{remark}

\def\be{\begin{equation}}
\def\ee{\end{equation}}

\def\ve{\varepsilon}

\begin{document}
\large

\title{On the Vinogradov mean value}
\date{\today}
\author{J.~Bourgain}
\address{J.~Bourgain, Institute for Advanced Study, Princeton, NJ 08540}
\email{bourgain@math.ias.edu}
\thanks{The author was partially supported by NSF grants DMS-1301619}

\begin{abstract}
A discussion of recent work of C.~Demeter, L.~Guth and the author of the
proof of the Vinogradov Main Conjecture using the decoupling theory for
curves.
\end{abstract}

\maketitle

\section
{Introduction and statements}

For $k, s\in\mathbb N$ and $x\in \mathbb R^k$, denote $\big(e(t)=e^{2\pi it}\big)$
$$
f_k(x, N)=\sum_{1\leq n\leq N} e(nx_1+n^2 x_2+\cdots+n^kx_k)\eqno{(1.1)}
$$
and
$$
J_{s, k}(N) =\int_{[0, 1]^k} |f_k(x, N)|^{2s} dx_1\cdots dx_k.\eqno{(1.2)}
$$
By orthogonality, $J_{s, k}(N)$ counts the number of integral solutions of the system
$$
n^j_1+\cdots + n^j_s=n^j_{s+1}+\cdots + n^j_{2s} \quad (1\leq j\leq k)\eqno{(1.3)}
$$
where $1\leq n_i\leq N$ \ $(1\leq i\leq 2s)$.

The evaluation of $J_{s, k}(N)$ is a central problem of importance to several classical issues in analytic number theory,
including the Waring problem, bounds on Weyl sums and zero-free regions for the Riemann zeta function.
The introduction of the mean value (1.2) and its significance to number theory go back to the seminal work of I.M.~Vinogradov (cf. \cite{Vi2}).
This approach is referred to as `Vinogradov's Method'.

Following T.~Wooley, we call `Main Conjecture' the statement
$$
J_{s, k}(N)\ll N^\ve(N^s+N^{2s-\frac 12k(k+1)}) \ \text { for all $\ve>0$}\eqno{(1.4)}
$$
and $2s =k(k+1)$ the critical exponent.
We note that indeed both $N^s$ and $N^{2s-\frac 12 k(k+1)}$ are obvious lower bounds (up to multiplicative constants).

Vinogradov's original argument \cite {Vi1} for estimating $J_{s, k}(N)$ was refined by means of Linnik's $p$-adic approach \cite{Li}
and the work of Karatsuba \cite {Ka} and Stechkin \cite{St}, leading to the following bound for $s\geq k$
$$
J_{s, k}(N)\leq  D(s, k) N^{2s-\frac 12 k(k+1)+\eta_{s, k}}\eqno{(1.5)}
$$
with
$$
\eta_{s, k}=\frac 12 k^2\Big(1-\frac 1k\Big)^{[s/k]} \text { and } \ D(s, k)=\min (k^{csk}, k^{ck^3}).\eqno{(1.6)}
$$
The latter leads to an asymptotic formula
$$
J_{s, k}(N)\sim C(s, k) N^{2s-\frac 12 k(k+1)}.\eqno{(1.7)}
$$
provided
$$
s\geq k^2 (2\log k+\log\log k+ 5)\eqno{(1.8)}
$$
(see \cite{A-C-K}).

Major progress towards the Main Conjecture was achieved by T.~Wooley based on his efficient congruencing method.

\begin
{theorem} \ (see Theorem 4.1 in \cite {W6}, based on \cite{W1}, \cite{W2}, \cite {F-W}, \cite {W3},\cite{W4}, \cite{W5}).

The Main Conjecture for $J_{s, k}(N)$ holds when 

(i) $k= 1, 2, 3$.

(ii) $1\leq s\leq D(k)$, where $D(4)=8, D(5)=10$, \ldots and
$$
D(k) =\frac 12 k(k+1)-\frac 13 k+O(k^{3/2})\eqno{(1.9)}
$$

(iii) $k\geq 3$ and $s\geq k(k-1)$
\end{theorem}

The reader is referred to the survey paper \cite{W6} for a detailed discussion.
It should be noted that prior to \cite {W5}, the Main Conjecture was only known for $k\leq 2$.

Based on  a more general harmonic analysis principle - the so-called `decoupling theorem' for curves - the full Main Conjecture was finally
established by C.~Demeter, L.~Guth and the author in the fall of 2015 (see \cite{BDG}).

\begin{theorem}
(\cite{BDG}).
The Main Conjecture for $J_{s, k}(N)$ holds.
\end{theorem}

For $s>\frac 12 k(k+1)$, the prefactor $N^\ve$ in (1.4) may be dropped and one has the asymptotic formula (1.7).
In what follows, we will review some consequences to the Waring problem and Weyl syms, improving on earlier results.
Next, we will formulate the underlying harmonic analysis result with a brief discussion (the reader will find complete proofs in \cite {BDG}) and
conclude with some further comments.

Concerning applications to the zeta-function, our work as it stands does not lead to further progress.
The reason for this is that we did not explore the effect of large $k$ (possibly depending on $N$) and in the present form is likely very poor.
A similar comment applies to Wooley's approach.

References in this paper are far from exhaustive and only serve the purpose of this expos\'e.

\section
{The Asymptotic Formula in Waring's Problem}

Denote $R_{s, k}(n)$ the number of representations of the positive integer $n$ as sum of $s$ $kth$ powers.
For $s$ sufficiently large, one has the asymptotic formula
$$
R_{s, k}(n) =\frac {\Gamma(1+\frac 1k)^s}{\Gamma(\frac sk)} {\frak S}_{s, k} (n) n^{\frac sk-1} + o(n^{\frac sk-1})\eqno{(2.1)}
$$
where
$$
{\frak S}_{s, k}(n) =\sum^\infty_{q=1} \sum^q_{\substack{a=1\\ (a, q)=1}} \Big(\frac 1q \sum^q_{r=1} e_q(ar^k)\Big) e_q(-na)\eqno{(2.2)}
$$
is the singular series.

Denote $\tilde G(k)$ the smallest integer $s$ for which (2.1) holds.
Based on heuristic applications of the circle method, one expects $\tilde G(k) =k+1$ for $k\geq 3$, but known results are far weaker.
The Vinogradov main value theorem  plays a crucial role in the minor arcs analysis (see in particular \cite {W7}).
Moreover, in \cite {W7} the implications of (1.4) to $\tilde G(k)$ are worked out, at that time conjectural.
Recording Theorem 4.1 in \cite {W7}, one obtains therefore the bound

\begin{theorem}
For $k\geq 3$,
$$
\tilde G(k) \leq k^2+1-\max_{\substack{1\leq j\leq k-1\\ 2^s\leq k^2}} \Big\lceil \frac{k_j-2^j}{k+1-j}\Big\rceil.\eqno{(2.3)}
$$
\end{theorem}
denoting $\lceil t\rceil $ the smallest integer no smaller then $t$.

In particular $\tilde G(4) =15, \tilde G(k) \leq k^2-2(k\geq 5), \tilde G(k) \leq k^2-3(k\geq 8), \ldots$
$$
\tilde G(k) \leq k^2+1-\Big[\frac {\log k}{\log 2}\Big] \qquad (k\geq 3).\eqno{(2.4)}
$$
This is an improvement of all previously known bounds on $\tilde G (k)$, except for Vaughan's $\tilde G(3)\leq 8$ (\cite{Vau1}).

As we will see later, the bound (2.4) may be further improved for large $k$, due to the fact that our results also enable a certain improvement in
Hua's lemma.

For the record, we note that Wooley obtained $\tilde G(5)\leq 28, \tilde G(6)\leq 43, \tilde G(7)\leq 61, \ldots$
$$
\tilde G(k) \leq \big(1.5407\ldots + o(1)\big) k^2 \ \text { for large $k$}.\eqno{(2.5)}
$$

\bigskip
\section
{Weyl Sums}

Recalling (1.1), Weyl's theorem states (see \cite{Vau2} for instance).

\begin{theorem}
(H.~Weyl). 
With the notation (1.1), assume $(a, q)=1$ and $\big| x_k-\frac aq\big|\leq \frac 1{q^2}$. Then
$$
|f_k(x, N)|\ll N^{1+\ve} (q^{-1}+ N^{-1}+qN^{-k})^{ 2^{1-k}}.\eqno{(3.1)}
$$
\end{theorem}

It is well-known that for large $k$, Vinogradov's method leads to substantially better results.
As a consequence of Theorem 2, one gets (cf. \cite{Vau2}).

\begin{theorem}
Again with the notation (1.1), let $k\geq 3$, $2\leq j\leq k$ and assume
$$
\Big|x_j-\frac aq\Big|\leq \frac 1{q^2}, (a, q)=1.
$$
Then
$$
|f_k(x, N)|\ll N^{1+\ve}(q^{-1}+N^{-1}+qN^{-j})^{\sigma(k)}\ \text { with } \sigma(k)=\frac 1{k(k-1)}.\eqno{(3.2)}
$$
\end{theorem}

Theorem 5 improves Weyl's bound and later refinements due to Heath-Brown \cite{H-B} and Robert-Sargos \cite {R-S} for $k\geq 7$.

Wooley had proven (3.2) with $\sigma(k) =\frac 1{2(k-1)(k-2)}$, see \cite {W6}.

\section
{The Decoupling Theorem for curves}

It turns out that in fact (1.4) is a consequence of a more general harmonic analysis principle that we discuss next.

Let $\Gamma =\{(t, t^2, \ldots, t^k\}:0\leq t\leq 1\}$ be the moment curve (or, more generally a non-degenerate curve in $\mathbb R^k$).
Given $g:[0, 1]\to \mathbb C$ and an interval $J\subset [0, 1]$, define the extension operator
$$
E_{_J} g(x)=\int_J g(t) e(tx_1+ t^2x_2+\cdots+ t^k x_k) dt \qquad x=(x_1, \ldots, x_k)\in\mathbb R^k.\eqno{(4.1)}
$$
Given a ball $B=B(c_B, R)$ in $\mathbb R^k$, denote $\omega_B$ the weight function
$$
\omega_B(x) =\Big(1+\frac {|x-c_B|}{R}\Big)^{-100k}.\eqno{(4.2)}
$$

\begin{theorem}
(\cite{BDG}).
Let $k\geq 2$ and $0<\delta\leq 1$.
For each ball $B\subset \mathbb R^k$ of radius at least $\delta^{-k}$, one has the inequality
$$
\Vert E_{[0, 1]} g\Vert_{L^{k(k+1)}(\omega_B)} \ll \delta^{-\ve} \Big(\sum_{J\subset [0, 1], |J|=\delta} \Vert E_J g\Vert^2_{L^{k(k+1)}(\omega_B)}\Big)^{\frac 12}\eqno{(4.3)}
$$
where $J$ runs over a partition of $[0, 1]$ in $\delta$-intervals.
\end{theorem}

\noindent
{\bf Remarks.}

(i) Decoupling inequalities of the type were previously established in \cite {BD1} for smooth hypersurfaces in $\mathbb R^k$ with non-vanishing
curvature.  In particular, the case $k=2$ of Theorem 6 already appears in \cite{BD1}.  We also refer the reader to \cite{BD1} for the analysis background of
the decoupling problem.

(ii) The exponent $k(k+1)$ in (4.3) is best possible.
Let us point out that there is a similar decoupling inequality for $2\leq p<k(k+1)$, though for $k\geq 3$ this is not just a consequence of interpolation.

(iii) Our decoupling inequalities for curves appear in \cite {B1}, \cite {BD2}, \cite {B2}.
In particular, the reader is referred to \cite{B2} for an application to exponential sums and the Lindel\"of hypothesis for the Riemann-zeta function.

(iv) The weight function $\omega_B$ (rather than $1_B$) is a (necessary) technical issue but will often be ignored in our later discussion for simplicity.

It is easy to deduce Theorem 2 from Theorem 6.  One first observes that the decoupling theorem implies the following discretized version.

\begin{theorem}
For each $1\leq n\leq N$, let $\frac {n-1}N< t_n<\frac nN$ and let $R>N^k$.
For each $p\geq 1$, one has
$$
\begin{aligned}
\Big\{\frac 1{|B_R|} &\int\Big|\sum^N_{n=1} a_n e(t_n x_1+t^2_n x_2+\cdots+ t_n^k x_k)\Big|^p \omega_{B_R}(x) dx_1\ldots dx_k\Big\}^{1/p}\ll\\
&N^\ve \Big(1+N^{\frac 12(1-\frac{k(k+1)}p)}\Big) \ \Big(\sum|a_n|^2\big)^{\frac 12}.
\end{aligned}
\eqno{(4.4)}
$$
\end{theorem}

\noindent
(the case $p<k(k+1)$ is obtained by interpolation with $p=2$ and $p> k(k+1)$ with the obvious $p=\infty$ bound).

Taking $a_n=1$ and $p=2s$, it follows from (4.4) that the system of inequalities
$$
|t^j_{n_1} +\cdots+ t^j_{n_s} - t^j_{n_{s+1}} -\cdots - t^j_{n_{2s}}|<N^{-k} \quad (1\leq j\leq k)\eqno{(4.5)}
$$
has at most $N^\ve(N^s+N^{2s-\frac {k(k+1)}{2}})$ solutions in $1\leq n_1, \ldots, n_{2s} \leq N$.
Specifying $t_n=\frac nN$, Theorem 2 follows immediately.

\section
{Elements of the Proof of Theorem 6}

Most techniques involved in proving decoupling theorems had previously been developed in the study of the restriction and Kakeya problems in harmonic analysis.
These include wave packet decomposition, parabolic rescaling and the use of multi-linear analysis.
In what follows, we make a few mostly superficial comments on how they appear in the context of curves.

\subsection
{Wave Packet Decomposition}

Let $J\subset [0, 1]$ be a small interval and $\tau =\{\gamma(t) =(t, t^2, \ldots, t^k); t\in J\}\subset\Gamma$ the corresponding arc.
Then, roughly speaking, $|E_{_J}g|$ may be viewed as `essentially constant' on translates of the geometric polar $\overset \circ \tau$ of the convex hull of $\tau$.
Thus if $|J|=\delta$, these are $\frac 1\delta \times\frac 1{\delta^2} \times\cdots \times \frac 1{\delta^k}$-boxes oriented according to the Frenet basis of $\Gamma$.

\subsection
{Parabolic Rescaling}

Take $k=2$ and $J=[t_0, t_0+\sigma]\subset [0, 1]$.
Write for $t= t_0+\sigma t' \in J$
$$
x_1 t+x_2 t^2 =x_1 t_0+x_2 t_0^2+\sigma(x_1+ 2x_2 t_0) t'+\sigma^2 x_2(t')^2\eqno{(5.1)}
$$
and make a change of variables $x_1' =\sigma(x_1+2t_0x_2), x_2'=\sigma^2 x_2$.

The map $(x_1, x_2)\mapsto (x_1', x_2')$ maps $B_R$ to an $\sigma R\times \sigma^2 R$ size ellipse which we cover with $\sigma^2R$-balls.

In general, $(x_1, \ldots, x_k)\mapsto (x_1', \ldots , x_k')$ maps $B_R$ to an ellipsoid covered by $\sigma^k R$-balls.
Next, denote $K_p(\delta)$ the best constant for which a decoupling inequality
$$
\Vert E_{[0, 1]} g\Vert_{L^p(B)}\leq K_p(\delta) \Big(\sum_{|J|=\delta} \Vert E_{_J}g\Vert^2_{L^p(B)}\Big)^{\frac 12}\eqno{(5.2)}
$$
with $B$ a $\delta^{-k}$-ball holds.
It follows then from the previous discussion that if $J\subset [0, 1], |J|=\sigma>\delta$, then (5.2) will hold with $K_p(\delta)$
replaced by $K_p(\frac\delta\sigma)$ if supp\,$g\subset J$.

\subsection
{Multilinear Analysis}

The reduction of (4.3), which in some sense is a linear statement, to multi-linear expressions is crucial as it allows us to exploit transversality.
This technique, which is basically simple, goes back to the joint work \cite{BG} of L.~Guth and the author.
All available results on decoupling make use of this procedure.

Continuing our high-level discussion, the left side of (4.3) will be replaced by certain multi-linear quantities which we describe next.
Define
$$
D_q(\Delta, \delta_1) =\prod^M_{i=1} \Big[\sum_{J\subset J_i, |J|=\delta_1} \Vert E_{_J} g\Vert^2_{L^q_{\#}(\Delta)} \Big]^{\frac 1{2M}}\eqno{(5.3)}
$$
where

$M=M_k$ is an appropriate integer $(M_2=2)$

$J_1, \ldots, J_M\subset [0, 1]$ are fixed $O(1)$-separated intervals

$\Delta=R$-ball, $R>\delta_1^{-1}$ and $L^q_{\#}(\Delta)$ is the normalized $L^q$-norm on $\Delta$.

Let $B$ be a (fixed) large ball and define further for $2\leq q\leq p$
$$
\tilde D_q (M, \delta_1) =\left[\begin{aligned} &\text {Average } \\
&\Delta = \text{$R$-ball} \subset B\end{aligned} \quad  D_q(\Delta, \delta_1)^p\right]^{\frac 1p}\eqno{(5.4)}
$$
Hence $\tilde D_p (R, \delta_1)\leq D_p(B, \delta_1)$.
The strategy is to bound $\tilde D_q(R, \delta_1)$ by gradually decreasing $\delta_1$ and increasing $R$.
Note that from the previous discussion, one has for $\delta<\delta_1, |B|> \delta^{-k}$
$$
D_p (B, \delta_1)\leq K_p \Big(\frac \delta{\delta_1}\Big) D_p(B, \delta).\eqno{(5.5)}
$$
Clearly, from basic orthogonality, if $\delta_1> \frac 1R$, then
$$
D_2(\Delta, \delta_1)\lesssim D_2\Big(\Delta, \frac 1R\Big)\eqno{(5.6)}
$$
\big(a rigorous justification requires in fact replacing $\Delta$ by a weight function $w_\Delta$ of the type (4.2)\big).

More generally, if $q\leq d(d+1), d<k$ and $R>\delta_1^{-d}$, one has
$$
D_q(\Delta, \delta_1) \ll R^\ve D_q(\Delta, R^{-\frac 1d})\eqno{(5.7)}
$$
by appealing to the decoupling theorem in dimension $d$ (exploiting only the variables $x_1, \ldots, x_k)$, assuming the latter already obtained.

We also note the following interpolation property, which is immediate from H\"older's inequality.
Let $q_1\leq q\leq q_2$ and $\frac 1q =\frac {1-\theta}{q_1} +\frac \theta{q_2}$.
Then
$$
D_q(\Delta, \delta_1)\leq D_{q_1} (\Delta, \delta_1)^{1-\theta} D_{q_2} (\Delta, \delta_1)^\theta\eqno{(5.8)}
$$
and similarly for $\tilde D_q$.

Next, the ball inflation, i.e. the increment of $R$, uses essentially transversality which comes with the multi-linear structure
of (5.3).
The main inequality writes
$$
\tilde D_{d\frac pk} (\delta_1^{-d}, \delta_1) \ll \tilde D_{d\frac pk} ( \delta_1^{-d-1}, \delta_1) \ \text{ for } \ 1\leq d< k\eqno{(5.9)}
$$
and follows from wave packet decomposition as explained in (4.1) and multi-linear Kakeya type estimates originating from the work \cite{BCT}.
Note that in (5.9) and keeping in mind (5.4), we are essentially trading an $L^p$-norm for an $L^{\frac dk p}$-norm.
This is possible by exploiting certain transversality properties.
The key result is the Brascamp-Lieb inequality that underlies the multi-linear Kakeya theory and we formulate next.

\begin{theorem}
(Brascamp-Lieb, see \cite{BDG} for related references).

Let $d\leq k$ and for $1\leq i\leq M$, let $V_i$ be a $d$-dimensional subspace of $\mathbb R^k$.
Denote $\pi_i:\mathbb R^k\to V_i$ the orthogonal projection.
We assume the following transversality condition
$$
\frac dk \dim V\leq \frac 1M \sum^M_{i=1} \dim (\pi_i V)\eqno{(5.10)}
$$
for all linear subspaces $V$ of $\mathbb R^k$.

Then the quantity
$$
\sup_{g_i \in L^1(V_i)} \ \frac {\Vert [ \prod^M_{i=1} |g_i \circ \pi_i|]^{\frac 1M}\Vert_{L^{k/d}(\mathbb R^k)}}{[\prod^M_{i=1} \Vert g_i\Vert_{L^1(V_i)}]^{\frac 1M}}
\eqno{(5.11)}
$$
is finite.
\end{theorem}

In the present application, the spaces $V_i$ are obtained as $V_i=[\gamma'(t_i), \ldots, \gamma^{(d)}(t_i)]$ with $t_i\in J_i$ $(1\leq i\leq M)$ introduced above and condition 
(5.10) for
appropriate $M$ results from the assumption that the curve $\Gamma$ is non-degenerate.

Let $p<k(k+1)$ be sufficiently close to $k(k+1)$.

Let $\delta_0= \delta^u$ with $u>0$ fixed and arbitrarily small.
Starting from $\tilde D_2 (\delta_0^{-1}, \delta_0)$, it follows from (5.9) that
$$
\tilde D_2(\delta_0^{-1}, \delta_0)\leq \tilde D_{\frac pk}(\delta_0^{-1}, \delta_0)\ll \tilde D_{\frac pk} (\delta _0^{-2}, \delta_0).
\eqno{(5.12)}
$$
Next, use (5.8) with $q=\frac pk, q_1=2, q_2=2\frac pk$ and (5.9) with $d=2$ to get
$$
\tilde D_{\frac pk}(\delta_0^{-2}, \delta_0)\ll \tilde D_2 (\delta_0^{-2}, \delta_0^2)^{1-\theta_1} \ \tilde D_{\frac{2p}k} (\delta_0^{-3}, \delta_0)^{\theta_1}\eqno{(5.13)}
$$
for some $0<\theta_1<1$.
The second factor in (5.13) is further processed interpolating between $q_1=6$ and $q_2=\frac {3p}k$,
Applying (5.7) with $d=2, q=6$ and (5.9) leads to
$$
\tilde D_{\frac {2p}{k}}(\delta_0^{-3}, \delta_0)\ll \tilde D_6(\delta_0^{-3}, \delta_0^{\frac 32})^{1-\theta_2} \ 
\tilde D_{\frac {3p}k} (\delta_0^{-4}, \delta_0)^{\theta_2}\eqno{(4.14)}
$$
for some $0<\theta_2< 1$. Next,
$$
\tilde D_6 (\delta_0^{-3}, \delta_0^{\frac 32})\leq \tilde D_2 (\delta_0^{-3}, \delta_0^3)^{1-\psi_2} \ \tilde D_{\frac {2p}k}(\delta_0^{-3\frac 32}, \delta_0^{\frac
32})^{\psi_2}\eqno{(5.15)}
$$
for some $0<\psi_2<1$.
The above are the first few steps of an interpolation scheme that together with inequality (5.5) and a bootstrap argument eventually permits us to estimate $K_p(\delta)\ll
\delta^{-\ve}$ for $p<k(k+1)$.

The sole purpose of the above discussion is to give the reader some sense of how the proof of Theorem 6 works, again referring to \cite {BDG} for the full account and further
references.

\section
{Some Further Comments}

\subsection
{An Improvement of Hua's Inequality}

We point out another arithmetical consequence of Theorem 6 related to Hua's lemma.
Recall the statement (of \cite {Vau2}).

\begin{theorem}
(Hua).  For $k\geq 1$, denote
$$
S(x)=\sum_{1\leq n\leq N} e(n^kx).\eqno{(6.1)}
$$
Then for $1\leq \ell\leq k$ we have
$$
\int_0^1 |S(x)|^{2^\ell} dx\ll N^{2^\ell-\ell+\ve} \ \text { for all } \ \ve>0.\eqno{(6.2)}
$$
\end{theorem}

We sketch the proof of the following

\begin{theorem}
Let $S(x)$ be defined by (6.1) and $s\leq k$ a positive integer.
Then
$$
\int_0^1 |S(x)|^{s(s+1)} dx \ll N^{s^2+\ve} \ \text { for all } \ \ve>0.\eqno{(6.3)}
$$
\end{theorem}

Clearly Theorem 10 improves upon Theorem 9 for $\ell \geq 5$.

\noindent
{\bf Proof of Theorem 10.}

We apply the decoupling theorem to the non-degenerate curve in $\mathbb R^s$
$$
\Gamma=\{(t^k, t^{s-1}, \ldots, t), 1\leq t\leq 2\}.\eqno{(6.4)}
$$
The discretized version analogous to Theorem 7 implies
$$
N^{-s} \int_{[-N, N]^s} \Big|\sum^{2N}_{n=N} e \Big(\Big(\frac nN\Big)^k x+\Big(\frac nN\Big)^{s-1} x_{s-1}+\cdots +\frac nN x_1\Big|^{s(s+1)}
dx_1\ldots dx_{s-1} dx \ll N^{\frac 12 s(s+1)+\ve}.\eqno{(6.5)}
$$
Rescaling and use of periodicity gives
$$ 
\int_{[-1, 1]}\int_{[0, 1]^{s-1}}\Big|\sum^{2N}_{n=N} e\Big(\frac {n^k}{N^{k-s}}x+n^{s-1} x_{s-1}+\cdots+ nx_1\Big)\Big|^{s(s+1)}dx_1\ldots dx_{s-1} dx\ll 
N^{\frac 12 s(s+1)+\ve}.\eqno{(6.6)}
$$
Denote $K_r =K_r(t)$ the kernel on $\mathbb T=\mathbb R/\mathbb Z$ which Fourier transform $\widehat K_r$ is
trapezoidal, satisfying $\widehat K_r(n)=1$ for $|n|\leq r$ and supp\,$\widehat K_r\subset [-2r, 2r]$.
Hence $\Vert K_r\Vert_1\leq 3$.
Multiply the integrand in (6.6) by
$$
K_{2N}(x_1) K_{2N^2}(x_2) \cdots K_{2N^{s-1}}(x_{s-1})\eqno{(6.7)}
$$
and perform the integration in $x_1, \ldots, x_{s-1}$.
Since (6.7) $\leq C_s\, N^{\frac 12 s(s-1)}$, it follows from (6.6) that
$$
\int_{[-1, 1]}\Big|\sum^{2N}_{n=N} e\Big(\frac {n^k}{N^{k-s}} x\Big)\Big|^{s(s+1)} dx \ll N^{s^2+\ve}.\eqno{(6.8)}
$$
Note that inequality (6.8) is essentially optimal and implies the weaker statement
$$
\int_0^1\Big|\sum^{2N}_{n=N} e(n^k x)\Big|^{s(s+1)} dx\ll N^{s^2+\ve}.\eqno{(6.9)}
$$
This proves (6.3).

Returning to the discussion in \S2 and \cite{W7}, we point out that in the treatment of the minor arcs in the circle method,
besides Vinogradov's inequality also Hua's lemma (Theorem 9) is involved in deriving Theorem 3 (see \S3  in \cite {W7}).
Hence Theorem 10 is expected to produce further improvements in bounding $\tilde G(k)$, which we discuss next (referring to \cite{W7} for details).

Following \cite {W7}, define the set $\mathcal M=\mathcal M_k$ of minor arcs as the set of real numbers $x\in [0, 1[$ with the property
that, whenever $a\in\mathbb Z, q\in\mathbb Z_+, (a, q)=1$ satisfy $|qx-a|\leq (2k)^{-1} N^{1-k}$, then $q>(2k)^{-1}N$.

Injecting (1.4) with $s=\frac 12 k(k+1)$ in Theorem 2.1 of \cite {W7} implies
$$
\int_{\mathcal M} |S(x)|^{k(k+1)} dx \ll N^{k^2-1-\ve}.\eqno{(6.10)}
$$
Inequality (6.10) is then interpolated with (6.2) or alternatively (6.3) in order to establish an inequality of the form
$$
\int_{\mathcal M} |S(x)|^{s_0} dx < N^{s_0-k-\tau}\eqno{(6.11)}
$$
for some $\tau>0$ and as small as possible exponent $s_0\in\mathbb Z_+$ that will provide a bound on $\tilde G(k)$.

Taking $s < k$ a parameter, let $s_0\in\mathbb Z_+, s(s+1)\leq s_0\leq k(k+1)$.
Interpolation between (6.10) and (6.3) gives
$$
\int_{\mathcal M}|S(x)|^{s_0} dx\ll N^{s_0-\eta+\ve}\eqno{(6.12)}
$$
with
$$
\eta =(1-a)(k+1)+as \ \text { and } \ a=\frac {k(k+1)-s_0}{k(k+1) -s(s+1)}.\eqno{(6.13)}
$$
Hence, in order to obtain (6.11), we are lead to the condition $\eta>k$, which translates in
$$
s_0>k^2- \frac {k-s-1}{k+1-s} s.\eqno{(6.14)}
$$
Consequently, we proved

\begin{theorem}
$$
\tilde G(k) \leq k^2+1-\max_{s\leq k} \Big\lceil s\frac{k-s-1}{k-s+1} \Big\rceil.\eqno{(6.15)}
$$
\end{theorem}

The reader will verify that (6.15) improves over Theorem 3 for $k>12$ and moreover implies that for large $k$
$$
\tilde G(k)< k^2-k +O(\sqrt k)\eqno{(6.16)}
$$
rather then (2.4).

\subsection
{Generalizations of Vinogradov's Inequality}
\hfil\break
Mean value estimates for multi-dimensional Weyl sums using efficient congruencing were obtained in \cite {PPW}.
One could reasonably expect that a complete understanding of decoupling phenomena for surfaces in $\mathbb R^k$ will
also lead to progress and perhaps optimal results in this more general setting.
Presently, we only reached a satisfactory understanding of decoupling for co-dimension one surfaces and for curves.
A decoupling theorem for 2-dimensional surfaces in $\mathbb R^k$ was established in \cite {BD3} implying in particular 
results on 2-dimensional cubic Weyl sums but that are likely not optimal.
The recent developments around curves obtained in \cite {BDG} almost surely will further contribute in this direction.


\begin{thebibliography}{XXXXX}

\bibitem[A-C-K]{A-C-K} Arkhipov, G.~I., Chubarikov, V.~N. and Karatsuba, A.~A. \, {\em Trigonometric sums in number theory and analysis}, Walker de Gruyter, Berlin, 2004.

\bibitem[B1]{B1} Bourgain, J.\, {\em Decoupling inequalities and some mean-value theorems}, to appear in Journal d'Analyse Mathematique, arXiv 1406.7862.

\bibitem [B2]{B2} Bourgain, J. \, {\em Decoupling, exponential sums and the Riemann zeta function}, to appear in JAMS, arXiv 1408.5794.

\bibitem[BCT]{BCT} Bennett, J., Carbery, A. and Tao, T. \, {\em On the multilinear restriction and Kakeya conjectures}, Acta Math. {\bf 196} (2006), no. 2, 261--302.


\bibitem[BD1]{BD1} Bourgain, J. and Demeter, C. \, {\em The proof of the $l^2$ Decoupling Conjecture}, Annals of Math. {\bf 182} (2015), no. 1, 351-389.

\bibitem[BD2]{BD2}  Bourgain, J. and Demeter, C. \, {\em Decouplings for curves and hypersurfaces with nonzero Gaussian curvature}, to appear in 
Journal d'Analyse Mathematique, arXiv 1409.1634.

\bibitem[BD3]{BD3} Bourgain, J. and Demeter, C. \, {\em Mean value estimates for  Weyl  sums in two dimensions},  arXiv 1509.05388.

\bibitem[BG] {BG}Bourgain, J. and Guth, L. \, {\em Bounds on oscillatory integral operators based on multilinear estimates}, Geom. Funct. Anal. {\bf 21} (2011), no. 6, 1239-1295.

\bibitem[BDG] {BDG} Bourgain, J., Demeter, C. and  Guth, L.\, {\em Proof of the main conjecture in Vinogradov's mean value theorem for
degrees higher then three}, arXiv 1512.01565 (2015).

\bibitem [F-W]{F-W} Ford, K.B. and  Wooley, T.D.\,
{\em On Vinogradov's mean value theorem: strongly diagonal behaviour via efficient congruencing}, Acta Math. {\bf 213} (2014), No 2, 199--236.

\bibitem[H-B]{H-B}  Heath-Brown, D.R. \, {\em Weyl's inequality, Hua's inequality, and Waring's problem}, J.~Londin Math. Soc. (20) {\bf 38} (1988), 396--414.

\bibitem[Ka]{Ka} Karatsuba, A.~A. \, {\em The mean value of the modulus of a trigonometric sum}, Izv. Akad.~Nauk SSSR Ser. Mat. {\bf 37} (1973), 1203--1227.

\bibitem[Li]{Li} Linnik, Yu.~V. \, {\em On Weyl's sums}, Mat. Sbornik (Rec. Math.) N.~S. {\bf 12} (1943), 28--39.

\bibitem[PPW]{PPW} Parsell, S.~T.,  Prendiville, S.~M. and Wooley, T.~D. \, {\em Near-optimal mean value estimates for multidimensional Weyl sums}, Geom.
Func. Anal. {\bf 23} (2013), no. 6, 1962--2024.

\bibitem[R-S]{R-S} Robert, O., Sargos, P. \, {\em Un th\'eor\`eme de moyenne pour les sommes d'exponentielles. Application \`a l'in\'egatite de Weil},
Publ.~Inst. Math. (Beograd) N.S. {\bf 67} (2000), 14--30.

\bibitem [St] {St} Stechkin, S. R. \, {\em On mean values of the modulus of a trigonometric sum}, Trudy Mat. Inst. Steklov {\bf 134} (1975), 283--309.

\bibitem [Vi1]{Vi1} Vinogradov, I.~M. \, {\em New estimates for Weyl sums}, Dokl. Acad. Nauk SSSR {\bf 8} (1935), 195--198.
 
\bibitem [Vi2]{Vi2} Vinogradov, I.~M. \, {\em The method of trigonometrical sums in the theory of numbers}, Trud. Inst. Mat. Steklov {\bf 23} (1947).
 
\bibitem [Vau1]{Vau1} Vaughan, R.~C.\, {\em On Waring's problem for cubes}, J. Reine Angew. Math. {\bf 365} (1986), 122--170.

\bibitem [Vau2]{Vau2} Vaughan, R.~C.\, {\em The Hardy-Littlewood method}, 2nd edition, Cambridge University Press, Cambridge, 1997.

\bibitem [W1]{W1} Wooley, T.~D. \, {\em Vinogradov's mean value theorem via efficient congruencing}, Ann. of Math. (2), {\bf 175} (2012), no 3, 1575--1627.

\bibitem[W2]{W2} Wooley, T. D.\, {\em   Vinogradov's mean value theorem via efficient congruencing, II}, Duke Math. J. {\bf 162} (2013), no 4, 673--730.

\bibitem [W3]{W3} Wooley, T.~D.\, {\em Multigrade efficient congruencing and Vinogradov's mean value theorem}, Proc. LMS (3), {\bf 111} (2015), no 3, 519--560.

\bibitem[W4]{W4} Wooley, T. D.\, {\em  Approximating the main conjecture in Vinogradov's Mean value Theorem}, arXiv 1401.2932.

\bibitem[W5]{W5} Wooley, T. D.\, {\em  The cubic case of the mean conjecture in Vinogradov's mean value theorem}, 
arXiv 1401.3150.

\bibitem [W6]{W6} Wooley, T. D.\, {\em Translation invariance, exponential sums and Waring's problem}, arXiv:1404.3508vi (2014).

\bibitem[W7]{W7} Wooley, T. D.\, {\em  The asymptotic formula in Waring's problem}, IMRN, Vol. 2012, No. 7, 1485-1504.




\end{thebibliography}
\end{document}